\documentclass[A4paper,11pt]{article}
\usepackage{latexsym}
\usepackage{mathrsfs}
\usepackage{amssymb}
\usepackage{amscd}
\usepackage[dvips]{graphicx}                  

\newtheorem{theorem}{Theorem}
\newtheorem{corollary}{Corollary}
\newtheorem{lemma}{Lemma}

\newenvironment{proof}[1][Proof]{\textbf{#1.} }{\ \rule{0.5em}{0.5em}}
\long\def\symbolfootnote[#1]#2{\begingroup%
	\def\thefootnote{$\;$}\footnote[#1]{$^*$#2}\endgroup}
\begin{document}
	
	\title{Ultrafilters without immediate predecessors in Rudin-Frolik order for regulars}
	\author{Joanna Jureczko\footnote{The author is partially supported by Wroc\l{}aw University of Science and Technology grant of K34W04D03 no. 8211104160.}}
\maketitle

\symbolfootnote[2]{Mathematics Subject Classification: Primary 03C25, 03E35, 03E55, 54E52,  54D35, 54E45, 54C10.

\hspace{0.2cm}
}

\begin{abstract}
	The aim of this paper is to construct ultrafilters without immediate predecessors in the Rudin-Frolik order in $\beta \kappa\setminus \kappa$, where $\kappa$ is a regular cardinal. This generalizes the problem posed by Peter Simon more than 40 years ago. 
\end{abstract}

\section{Introduction}

The Rudin-Frolik ordering of ultrafilters has a long history. It was defined by Z. Frolik in \cite{ZF} who used it to prove that $\beta\omega\setminus \omega$ is not homogeneous. M.E. Rudin, who nearly defined this ordering in cite \cite{MR1}, is credited with being the first to notice that the relationship between filters she used is actually ordering. D. Booth in \cite{DB} showed that this relation is a partial ordering of the equivalence classes, that is a tree, and that it is not well-founded.

The author defined and studied the partial orders on the type of points in $\beta\omega$ and in $\beta\omega \setminus \omega$ in \cite{MR}. These definitions were used later in \cite{BB} and \cite{BE2}.

E. Butkovi\v cov\'a between 1981 and 1990, published a number of papers concerning ultrafilters in the Rudin-Frolik order in $\beta \omega\setminus \omega$. In \cite{BB} with L. Bukovsk\'y and \cite{BE2} she constructed an ultrafilter on $\omega$ with a countable set of its predecessors. In \cite{BE1} she constructed ultrafilters without immediate predecessors. In \cite{BE3} Butkovi\v cov\'a showed that there exists in Rudin-Frolik order an unbounded chain order-isomorphic in $\omega_1$. In \cite{BE4} she proved that there is a set of $2^{2^{\aleph_0}}$ ultrafilters incomparable in Rudin-Frolik order which is bounded from below and no its subset of cardinality more than one has an infimum. In \cite{BE5} Butkovi\v cov\'a proved that for every cardinal between $\omega$ and $\mathfrak{c}$ there is a strictly decreasing chain without a lower bound. In most of these papers, the method presented in \cite{KK}.

In 1976, A. Kanamori published a paper \cite{AK} in which, among other things, it showed that the Rudin-Frolik tree cannot be very high if one considers it over a measurable cardinal. Moreover, in the same paper he left a number of open problems regarding the Rudin-Frolik order. Recently, M. Gitik in \cite{MG} answered some of them by using meta-mathematical methods. The solutions to some of the problems from \cite{AK} presented in combinatorial methods are in preparation, (\cite{JJ_kanamori}).

Since the Rudin-Frolik order is still not well researched, especially for ultrafilters over cardinals greater than $\omega$, it seems worthwhile to look at them for such cardinals. For this purpose, special techniques are needed. 
 
In 2001, Baker and Kunen published in \cite{BK} a very useful method which can be recognized as a generalization of the method presented in \cite{KK}. It is worth emphasizing that both methods, (from \cite{KK} and \cite{BK}), provide useful "technology" for keeping the transfinite construction of an ultrafilter not finished before $\mathfrak{c}$ steps, (see \cite{KK}), and $2^\kappa, $ for $\kappa$ being infinite cardinal, (see \cite{BK}), but the second method has some limitations, among others $\kappa$ must be regular. As it turned out, the method from \cite{BK} can be useful in keeping the results for the Rudin-Frolik order but for $\kappa$-regular. Due to the lack of adequate useful method for $\kappa$-singular, the similar results but for singulars are still open questions. Based on \cite{BK} and papers by Butkovi\v cov\'a, the investigations into the Rudin-Frolik order for regulars have further consequences, which will be presented in \cite{JJ_order1, JJ_order2}. \\ 
The starting point of this paper, based on \cite{BE1}, is the following problem posed by Peter Simon, known as Simon's problem: \\\\ \textit{Does there exist a non-minimal ultrafilter in the Rudin-Frolik order of $\beta \omega \setminus \omega$ without immediate predecessor?} \\\\ One does not know how old the question is, but it was mentioned in \cite{BB} in 1981 and solved in \cite{BE1} in 1982. As was mentioned above, the main method used in \cite{BE1} is the method introduced by K. Kunen in \cite{KK} for keeping the transfinite construction of an ultrafilter from finishing before continuum steps.

As far as one knows, nobody considers the Simon's problem in the general case, i.e. for $\beta \kappa$, for arbitrary cardinal $\kappa$. This is the main idea for checking whether such a generalization is possible. The answer is yes, but for $\kappa$-regular. It is dictated by an adequate method based on the one presented in \cite{BK}, which just works for regulars.

The paper is organized as follows: in Section 2, there are presented definitions and previous facts needed for proving the results in further parts of this paper. Section 3 contains main results preceded by auxiliary results. The paper is finished with the open problem for singulars.

\section{Definitions and previous results}

\textbf{2.1.} In the whole paper, we assume that $\kappa$ is an infinite cardinal. Then $\beta\kappa$ means the \v Cech-Stone compactification, where $\kappa$ has the discrete topology. Hence, $\beta\kappa$ is the space of ultrafilters on $\kappa$ and $\beta \kappa \setminus \kappa$ is the space of non-principal ultrafilters on $\kappa$. We will write $u(\kappa)$ to denote the space of uniform ultrafilters, i.e. $\mathcal{F} \in u(\kappa)$ iff $|A| = \kappa$ for all $A \in \mathcal{F}$. \\\\ \textbf{2.2.} A set $\{\mathcal{F}_\alpha \colon \alpha < \kappa\}$ of filters on $\kappa$ is \textit{$\kappa$-discrete} iff there is a partition $\{A_\alpha \colon \alpha < \kappa\}$ of $\kappa$ such that $A_\alpha \in \mathcal{F}_\alpha$ for each $\alpha < \kappa$. \\The following fact, taken from \cite{MR}, is easily generalized \\\\ \textbf{Fact 1.} Let $\kappa$ be a regular infinite cardinal. If $X, Y$ are $\kappa$-discrete sets of ultrafilters on $\kappa$ and $\mathcal{G} \in \overline{X} \cap \overline{Y}$, then $$\mathcal{G} \in \overline{X \cap Y} \cup (\overline{X \cap (\overline{Y} \setminus Y)}) \cup \overline{(Y \cap (\overline{X}\setminus X))}.$$ \textbf{2.3.} Let $\mathcal{F}, \mathcal{G}\in \beta \kappa$. We define \textit{Rudin-Frolik order} as follows $$\mathcal{F} \leqslant_{RF} \mathcal{G} \textrm{ iff } \mathcal{G} = \Sigma(X, \mathcal{F})$$ for some $\kappa$-discrete set $X = \{\mathcal{F}_\alpha \colon \alpha < \kappa\} \subseteq \beta\kappa,$ where $$\Sigma(X, \mathcal{F}) = \{A \subseteq \kappa \colon \{\alpha < \kappa \colon A \in \mathcal{F}_\alpha\}\in \mathcal{F}\}.$$ We define $$\mathcal{F} =_{RF} \mathcal{G} \textrm{ iff } \mathcal{F} \leqslant_{RF} \mathcal{G} \textrm{ and } \mathcal{G} \leqslant_{RF} \mathcal{F}$$ $$\mathcal{F} <_{RF} \mathcal{G} \textrm{ iff } \mathcal{F} \leqslant_{RF} \mathcal{G} \textrm{ and } \mathcal{F} \not =_{RF} \mathcal{G}.$$ \textbf{2.4.} An ultrafilter $\mathcal{G} \in \beta \kappa \setminus \kappa$ is \textit{non-minimal} in Rudin-Frolik order iff there exists a $\kappa$-discrete set $X \subseteq \beta \kappa \setminus \kappa$ of ultrafilters such that $\mathcal{G} \in \overline{X} \setminus X$. 

An ultrafilter $\mathcal{G} \in \beta \kappa \setminus \kappa$ \textit{has an immediate predecessor} $\mathcal{F}$ in Rudin-Frolik order iff $\mathcal{G} = \Sigma(X, \mathcal{F})$ for some $\kappa$-discrete set $X$ of minimal ultrafilters in Rudin-Frolik order. 

A \textit{Simon point} in $\beta \kappa \setminus \kappa$ is an ultrafilter $\mathcal{G} \in \beta \kappa \setminus \kappa$ which is non-minimal and has no immediate predecessors. \\\\ \textbf{2.5.}Let $\tau$ and $\kappa$ be infinite cardinals. A set of filters $\{\mathcal{F}_{\xi, \zeta} \colon \xi < \tau, \zeta < \kappa\}$ is \textit{stratified} iff \begin{itemize} \item [(1)] $\{\mathcal{F}_{\xi, \zeta} \colon \zeta < \kappa\}$ is $\kappa$-discrete for each $\xi < \tau$, \item [(2)] for all $A \in \mathcal{F}_{\xi, \zeta}$, $\xi < \tau, \zeta < \kappa$ and each $\nu$ such that $\xi < \nu < \tau$ $$|\{\mu < \kappa \colon A \in \mathcal{F}_{\nu, \mu}\}| = \kappa.$$ \end{itemize} 

\textbf{2.6.} Let $X =\{\mathcal{F}_{\xi, \zeta} \colon \xi < \tau, \zeta < \kappa\}$ be a stratified set of filters and let $W$ be its subset. We define \begin{itemize} \item [(1)] $W(0) = W$, \item [(2)] $W(\gamma) = \bigcup_{\beta < \gamma} W(\beta)$, for $\gamma$ - limit, \item [(3)] $W(\gamma+1) = W(\gamma) \cup \{\mathcal{F}_{\xi, \zeta} \colon \exists_{\eta> \gamma}\ \exists_{{A \in \mathcal{F}_{\xi, \zeta}}} \{\mathcal{F}_{\eta, \nu}\colon A \in \mathcal{F}_{\eta, \nu}\}\subseteq W(\gamma)\}$, \item [(4)] $\tilde{W} = \bigcup_{\gamma < 2^\kappa} W(\gamma)$. \end{itemize} 

Intuitively, the above construction is used to select only certain filters from $X$ with the desired property (see e.g. property (P) in 2.8), and then add (inductively) to the set $ W $ only those filters outside $ W $ which satisfy the condition (2) in the definition in 2.5. This construction will be used to define property (P), (see 2.8), the formulation of which would not be possible by taking the entire ultrafiter family into account. \\\\ \textbf{2.7.} A function $\hat\varphi \colon [\kappa^+]^{<\omega} \to [\kappa]^{< \omega}$ is \textit{$\kappa$-shrinking} iff \begin{itemize} \item [(1)] $p\subseteq q$ implies $\hat{\varphi}(p) \subseteq \hat{\varphi}(q)$, for any $p, q \in [\kappa^+]^{<\omega}$, \item [(2)] $\hat{\varphi}(0) = 0$. \end{itemize} 

A \textit{step family} (over $\kappa$, with respect to $\hat{\varphi}$) is a family of subsets of $\kappa$, $$\{E_t \colon t \in [\kappa]^{<\omega}\} \cup \{A_\alpha \colon \alpha < \kappa^+\}$$ satisfying the following conditions: \begin{itemize} \item [(1)] $E_s\cap E_t = \emptyset$ for all $s, t \in [\kappa]^{<\omega}$ with $s \not =t$, \item [(2)] $|\bigcap_{\alpha \in p} A_\alpha \cap \bigcup_{t\not \supseteq \hat{\varphi}(p)}E_t| < \kappa$ for each $p \in [\kappa^+]^{<\omega}$, \item [(3)] if $\hat{\varphi}(p) \subseteq t$, then $|\bigcap_{\alpha\in p} A_\alpha \cap E_t|=\kappa$ for each $p \in [\kappa^+]^{<\omega}$ and $t \in [\kappa]^{<\omega}$. \end{itemize} 

Let $I$ be an index set and $\mathcal{F}$ be a filter on $\kappa$. The family $$\{E_t^i \colon t \in [\kappa]^{<\omega}, i \in I\} \cup \{A_\alpha^i \colon \alpha < \kappa^+, i \in I\}$$ is an \textit{independent matrix of $|I|$ step-families} (over $\kappa$) with respect to (shortly w.r.t.) $\mathcal{F}, \hat{\varphi}$ iff \begin{itemize} \item [(1)] for each fixed $i \in I$, $\{E_t^i \colon t \in [\kappa]^{<\omega}\} \cup \{A_\alpha^i \colon \alpha < \kappa^+\}$ is a step-family, \item [(2)] if $n \in \omega, p_0, p_1, ..., p_{n-1} \in [\kappa^+]^{<\omega}, t_0, t_1, ..., t_{n-1} \in [\kappa]^{<\omega}$, $i_0, i_1, ..., i_{n-1} \in I$ with $i_k\not = i_m, k\not = m$ and $\hat{\varphi}(p_k) \subseteq t_k$, then $$\kappa\setminus \big[ \bigcap_{k=1}^{n-1}(\bigcap_{\alpha\in p_k}A^{i_k}_{\alpha} \cap E^{i_k}_{t_k})\big] \not\in \mathcal{F}.$$ \end{itemize} \noindent \textbf{Fact 2 (\cite{BK}).} If $\kappa$ is a regular cardinal and $\hat{\varphi}$ is a $\kappa$-shrinking function, then there exists an independent matrix of $2^\kappa$ step-families over $\kappa$ with respect to the filter $\{D \subseteq \kappa \colon |\kappa \setminus D| < \kappa\}$, $\hat{\varphi}$. \\\\ \textbf{2.8.} Let $\{D_t \colon t < [\kappa]^{<\omega}\}$ be a partition of $\kappa$ and let $\hat{\varphi}$ be a $\kappa$-shrinking function. A stratified set of filters $\{\mathcal{F}_{\xi, \zeta} \colon \xi < \tau, \zeta < \kappa\}$ \textit{satisfies property (P)} iff $\mathcal{F}_{\nu, \mu} \not \in \tilde{W}$ implies $$\exists_{\{W_\eta \colon \eta < \kappa^+\} \subseteq \mathcal{F}_{\nu, \mu}} \forall_{p \in [\kappa^+]^{<\omega}}\ |\bigcap_{\eta \in p} W_\eta \cap D_{\hat{\varphi}(p)}| < \kappa,$$ where $$W = \{\mathcal{F}_{\xi, \zeta} \colon \exists_{t < [\kappa]^{<\omega}} D_t \in \mathcal{F}_{\xi, \zeta}\}.$$

\section{Auxiliary Results}

	\begin{lemma} Let $\kappa$ be a regular cardinal and let $\hat{\varphi}$ be a $\kappa$-shrinking function. Then there exists a stratified set of ultrafilters $\{\mathcal{F}_{\xi, \zeta} \colon \xi, \zeta < \kappa\}$ on $\kappa$ which satisfies property (P) for each partition of $\kappa$. \end{lemma}
 
\begin{proof} Accept the following notation: \begin{itemize} \item $\mathcal{FR}(\kappa) = \{A \subset \kappa \colon |\kappa\setminus A|< \kappa\},$ \item $[A, B, C,...]$ means a filter generated by $A, B, C, ...$, \item $A\subseteq^*B$ means $|A \setminus B|< \kappa$. \end{itemize} By Fact 2, fix a matrix $$\{E_t^i \colon t \in [\kappa]^{<\omega}, i \in 2^\kappa\} \cup \{A_\alpha^i \colon \alpha < \kappa^+, i \in 2^\kappa\}$$ of $2^\kappa$ step-families over $\kappa$ independent with respect to the filter $\mathcal{FR}(\kappa)$. Assume that for each $i \in 2^\kappa$ we have \begin{itemize} \item [(a)] $E^i_s\cap E^i_t = \emptyset$ for all $s, t \in [\kappa]^{<\omega}$ with $s \not =t$, \item [(b)] $|\bigcap_{\eta \in p} A^i_\eta \cap \bigcup_{t\not \supseteq \hat{\varphi}(p)}E^i_t| < \kappa$ for each $p \in [\kappa^+]^{<\omega}$, \item [(c)] $|\bigcap_{\eta\in p} A^i_\eta \cap E^i_t|=\kappa$ for each $p \in [\kappa^+]^{<\omega}$ and $t \in [\kappa]^{<\omega}$ with $\hat{\varphi}(p) \subseteq t$, \item [(d)] $\bigcup\{E^i_t \colon t \in [\kappa]^{< \omega}\} = \kappa,$ \item [(e)] $|\bigcap_{\eta \in p} A^i_\eta \setminus \bigcup_{t \supseteq \hat{\varphi}(p)}E^i_t| < \kappa$ for each $p \in [\kappa^+]^{<\omega}.$ \end{itemize} The conditions $(a)-(c)$ are taken from the definition of step families, but $(d)-(e)$ can be easily obtained by "shrinking" $A^i_\alpha$ to $A^i_\alpha \subseteq \bigcup \{E_t^i \colon t \in [\kappa]^{<\omega}\}$ and then expanding $E^i_t$ so that $\{E^i_t \colon t \in [\kappa]^{< \omega}\}$ is a partition of $\kappa$.

Note that $(b)$ and $(d)$ imply $(e)$. Moreover, the condition $(b)$ is still preserved after expanding $\{E_t^i \colon t \in [\kappa]^{<\omega}\}$ to a partition of $\kappa$. \\Indeed. If there are $p_0 \in [\kappa^+]^{<\omega}$ and $i_0 \in 2^\kappa$ such that $$|\bigcap_{\eta \in p_0} A_\eta^{i_0} \cap \bigcup_{t\not \supseteq \hat{\varphi}(p_0)}E_t^{i_0}| = \kappa,$$ then $|\bigcup_{t\not \supseteq \hat{\varphi}(p_0)}E_t^{i_0}| = \kappa.$ Then, by $(d)$ and $(a)$, there would exist $t_0 \supseteq \hat{\varphi}(p_0)$ such that $|E_{t_0}^{i_0}| < \kappa$. Hence $$|\bigcap_{\eta\in p_0} A^{i_0}_\eta \cap E^{i_0}_{t_0}|<\kappa.$$ which contradicts $(c)$. 
Let $\{Z_\alpha \colon \alpha < 2^\kappa\}$ be the family of all subsets of $\kappa$. Let $\mathcal{D} = \{\mathcal{D}^\alpha \colon \alpha < 2^\kappa\}$ be a set of all partitions of $\kappa$, where $\mathcal{D}^\alpha = \{D^\alpha_s \colon s < [\kappa]^{\omega}\}$ can occur $2^\kappa$ many times in this enumeration. 

Now, we will construct $\mathcal{F}^0_{\xi, \zeta}$ and $I_0$. For this purpose, fix a partition $$G = \{G_\gamma \colon \gamma< \kappa, |G_\gamma| = \kappa\}$$ of $\kappa$ and define a family $\{B_{\xi, \zeta} \colon \zeta < \kappa\}$ of pairwise disjoint sets in the sense $B_{\xi, \zeta}\cap B_{\xi, \zeta'} = \emptyset$, whenever $\zeta \not = \zeta'$ 
For any $\zeta < \kappa$ we will define $B_{0,\zeta}$ as follows. Let $p \in [\kappa^+]^{<\omega}$ be such that $\zeta \in p$ and denote such $p$ by $p_\zeta$. Fix $t \in [\kappa]^{<\omega}$ such that $t \supseteq \hat{\varphi}(p_\zeta)$. Such $t$ is denoted by $t_\zeta$. Take $$B_{0,\zeta} = (\bigcap_{\eta \in p_\zeta}A^0_\eta \cap E^0_{t_\zeta}) \setminus \bigcup_{\mu < \zeta}(\bigcap_{\eta \in p_\mu }A^0_\eta\cap E^0_{t_\mu}).$$ By the above construction, $B_{0,\zeta}, \zeta < \kappa$ are pairwise disjoint and well defined, (which follows from $(a)-(e)$).

Assume that we have constructed $B_{\delta, \zeta}$ for some $\delta< \xi< \kappa$. \\ If $\xi$ is limit, then set $B_{\xi, \zeta} = \bigcap_{\delta < \xi} B_{\delta, \zeta}$. \\ Now we will construct $B_{\xi, \zeta}$ for $\xi = \delta+1$. Observe that since $G$ is a partition of $\kappa$ then $\zeta \in G_\gamma$ for some $\gamma < \kappa$. Then take $$B_{\xi,\zeta} = B_{\delta, \gamma}\cap [(\bigcap_{\eta \in p_\zeta}A^\xi_\eta \cap E^\xi_{t_\zeta}) \setminus \bigcup_{\mu < \zeta}(\bigcap_{\eta \in p_\mu }A^\xi_\eta\cap E^\xi_{t_\mu})].$$ Thus, for any $\xi < \kappa$ we have constructed the family $\{B_{\xi, \zeta} \colon \zeta < \kappa\}$ of the required property. 
Now, for any $\xi, \zeta < \kappa$ define $$\mathcal{F}^0_{\xi, \zeta} =[\mathcal{FR}(\kappa), \{B_{\xi, \zeta}\}, \{\kappa \setminus B_{\nu, \eta} \colon \nu > \xi, \eta < \kappa\}]$$ and set $I_0 = 2^\kappa \setminus \kappa$. It is easy observation that $$\{\mathcal{F}^0_{\xi, \zeta} \colon \xi, \zeta < \kappa\}$$ is a stratified set of filters. Obviously, $$\{E_t^i \colon t \in [\kappa]^{<\omega}, i \in I_0\} \cup \{A_\alpha^i \colon \alpha < \kappa^+, i \in I_0\}$$ of step-families is independent w.r.t. $\mathcal{F}^0_{\xi, \zeta}, \hat\varphi.$

Now, we will construct filters $\mathcal{F}^\alpha_{\xi, \zeta}$ and indexed sets $I_\alpha$ by induction on $\alpha < 2^\kappa$ steps fulfilling the properties \begin{itemize} \item [(1)] $\mathcal{F}^0_{\xi, \zeta}$ and $I_0$ as are done above, \item [(2)] $\mathcal{F}^\alpha_{\xi, \zeta}$ is a filter on $\kappa$, $I_\alpha\subset 2^\kappa$ and the matrix $$\{E_t^i \colon t \in [\kappa]^{<\omega}, i \in I_\alpha\} \cup \{A_\alpha^i \colon \alpha < \kappa^+, i \in I_\alpha\}$$ of remaining step-families is independent w.r.t. $\mathcal{F}^\alpha_{\xi, \zeta}, \hat{\varphi}$, \item [(3)] $\mathcal{F}^\alpha_{\xi, \zeta} = \bigcup_{\beta < \alpha} \mathcal{F}^\beta_{\xi,\zeta}$, $I_\alpha = \bigcap_{\beta<\alpha} I_\beta$, for $\alpha$ - limit, \item [(4)] $\mathcal{F}^\beta_{\xi, \zeta} \subseteq \mathcal{F}^\alpha_{\xi, \zeta}$, $I_\beta \supseteq I_\alpha$, whenever $\beta < \alpha$, \item [(5)] $I_\alpha \setminus I_{\alpha + 1}$ is finite, \item [(6)] if $\alpha \equiv 0\ (mod\ 2)$, then either $Z_\alpha \in \mathcal{F}^\alpha_{\xi,\zeta}$ or $\kappa \setminus Z_\alpha \in \mathcal{F}^\alpha_{\xi, \zeta}$, \item [(7)] if $\alpha \equiv 1\ (mod\ 2)$, then the set $\{\mathcal{F}^\alpha_{\xi, \zeta}\colon \xi, \zeta < \kappa\}$ of filters is stratified and fulfills property $(P)$. \end{itemize} 
Then, take $$\mathcal{F}_{\xi, \zeta} = \bigcup_{\alpha < 2^\kappa} \mathcal{F}^\alpha_{\xi, \zeta}.$$ Thus, $\{\mathcal{F}_{\xi, \zeta} \colon \xi, \zeta < \kappa\}$ will be the required set of ultrafilters.

Assume that $\mathcal{F}^\beta_{\xi, \zeta}$ and $I_\beta$ have been constructed for some $\beta <\alpha < 2^\kappa$. The limit step is done. We show how to obtain $\mathcal{F}^{\alpha}_{\xi, \zeta}$ and $I_{\alpha}$ for $\alpha = \beta+1$. \\ \\ \underline{Case $\alpha\equiv 0\ (mod \ 2)$}. 
If $\mathcal{R} = [\mathcal{F}^\beta_{\xi, \zeta}, Z_\beta]$ is a proper filter and the matrix of step-families $$\{E_t^i \colon t \in [\kappa]^{<\omega}, i \in I_\beta\} \cup \{A_\eta^i \colon \eta < \kappa^+, i \in I_\beta\}$$ is independent w.r.t. $\mathcal{R}, \hat{\varphi}$, then put $\mathcal{F}^\alpha_{\xi, \zeta}= [\mathcal{F}^\beta_{\xi, \zeta}, Z_\beta]$ and $I_\alpha = I_\beta$.

Otherwise, fix $n \in \omega$, distinct $i_k \in I_\alpha$ and $\hat{\varphi}(p_k)\subseteq t_k,$ for $ k < n$, such that $$\kappa \setminus [Z_\alpha \cap \bigcap_{k=0}^{n-1}(\bigcap_{\eta\in p_k}A^{i_k}_{\eta} \cap E^{i_k}_{t_k})]\in \mathcal{F}^\alpha_{\xi, \zeta}.$$ Then, put $$\mathcal{F}^{\alpha}_{\xi, \zeta} = [\mathcal{F}^{\beta}_{\xi, \zeta}, \{A^{i_k}_{p_k} \colon 0 \leqslant k \leqslant n-1\}, \{E^{i_k}_{t_k} \colon 0 \leqslant k \leqslant n-1\}]$$ $$I_{\alpha} = I_{\beta}\setminus \{i_k \colon 0 \leqslant k \leqslant n-1\}.$$ Then $\kappa \setminus Z_\alpha \in \mathcal{F}^{\alpha}_{\xi, \zeta}.$

To show that $(2)$ is fulfilled for $\mathcal{F}^{\alpha}_{\xi, \zeta}$ and $I_\alpha$ it is enough to observe that each element of $\mathcal{F}^{\alpha}_{\xi, \zeta}$ is of the form $$A \cap (\bigcap_{\eta \in p_k}A^{i_k}_\eta \cap E^{i_k}_{t_k})$$ for some $A \in \mathcal{F}^{\alpha}_{\xi, \zeta} $ and $ 0 \leqslant k \leqslant n-1$. Then, $$\bigcap_{\eta \in p_k}A^{i_k}_\eta \cap E^{i_k}_{t_k} \supseteq\bigcap_{k=0}^{n-1}(\bigcap_{\eta\in p_k}A^{i_k}_{\eta} \cap E^{i_k}_{t_k}) $$ and $$\bigcap_{k=0}^{n-1}(\bigcap_{\eta\in p_k}A^{i_k}_{\eta} \cap E^{i_k}_{t_k}) \in \mathcal{F}^{\beta}_{\xi, \zeta}.$$ Thus, $(2)$ for $\mathcal{F}^{\alpha}_{\xi, \zeta}$ and $I_\alpha$ is fulfilled which follows from $(2)$ for $\mathcal{F}^{\beta}_{\xi, \zeta}$ and $I_\beta$. \\ \\ \underline{Case $\alpha\equiv 1\ (mod \ 2)$}.

If $D^\beta_s \in \mathcal{F}^\beta_{\xi, \zeta}$ for some $s \in [\kappa]^{<\omega}$, then put $\mathcal{F}^\alpha_{\xi, \zeta} = \mathcal{F}^\beta_{\xi, \zeta}$ and $Z_\alpha = I_\beta$. 
Otherwise, we have (by (6)) $\kappa\setminus D^\beta_s \in \mathcal{F}^\beta_{\xi, \zeta}$ for each $s \in [\kappa]^{<\omega}$. 

Consider sequences $\langle V^\beta_s \colon s \in [\kappa]^{<\omega}\rangle$ such that $V^\beta_t \subseteq V^\beta_s$ whenever $t \subseteq s$ and $D^\alpha_s \subseteq V^\alpha_s$ for any $s \in [\kappa]^{<\omega}$. 

Now, choose $i \in I_\beta$ such that $$\bigcup\{E^i_t\setminus V^i_s \colon s, t \in [\kappa]^{<\omega}, t \subseteq s\} \not = \emptyset.$$ Set $$W^\beta_\delta = A^i_\delta \cap \bigcup\{E^i_t\setminus V^i_s \colon s, t \in [\kappa]^{<\omega}, t \subseteq s\}.$$ Then, set $I_\alpha = I_\beta \setminus\{i\}$ and $\mathcal{F}^\alpha_{\xi, \zeta} = [\mathcal{F}^\beta_{\xi, \zeta}, \{W^\beta_\delta \colon \delta < \kappa^+\}].$

To show that $(2)$ holds for $\mathcal{F}^\alpha_{\xi, \zeta}$ and $I_\alpha$ it is enough to observe that each element of $\mathcal{F}^\alpha_{\xi, \zeta}$ is of the form $$A \cap \bigcap_{\delta \in p}W^\alpha_\delta$$ for some $A\in \mathcal{F}^\alpha_{\xi, \zeta}$ and $p \in [\kappa^+]^{<\omega}$. But $$\bigcap_{\delta\in p} W^\alpha_\delta = \bigcap_{\delta\in p} A^i_\delta \cap \bigcup\{E^i_t\setminus V^i_s \colon s, t \in [\kappa]^{<\omega}, t \subseteq s\} \cap (\kappa\setminus D^\alpha_{\hat{\varphi}(p)})$$ and $\kappa\setminus D^\alpha_{\hat{\varphi}(p)} \in \mathcal{F}^\beta_{\xi, \zeta}$. Thus, $(2)$ for $\mathcal{F}^{\alpha}_{\xi, \zeta}$ and $I_\alpha$ is fulfilled which follows from $(2)$ for $\mathcal{F}^{\beta}_{\xi, \zeta}$ and $I_\beta$. 

Using a similar argument, it is easy to check that $\{\mathcal{F}^\alpha_{\xi, \zeta} \colon \xi, \zeta < \kappa\}$ is a stratified set of filters. 

Now, we show that $\{\mathcal{F}^\alpha_{\xi, \zeta} \colon \xi, \zeta < \kappa\}$ fulfills property $(P)$. Indeed, for each $p \in [\kappa^+]^{<\omega}$ we have $$\bigcap_{\delta \in p} W^\alpha_\delta = \bigcap_{\delta\in p} A^i_\delta \cap \bigcup\{E^i_t\setminus V^i_s \colon s, t \in [\kappa]^{<\omega}, t \subseteq s\} $$ (by $(e)$) $$\subseteq^* \bigcup_{t \supseteq\hat{\varphi}(p)}E^i_t \cap \bigcup\{E^i_t\setminus V^i_s \colon s, t \in [\kappa]^{<\omega}, t \subseteq s\} $$ (by $(a)$) $$= \bigcup_{t \supseteq \hat{\varphi}(p)} E^i_t \setminus V^\alpha_t \subseteq \bigcup_{t \supseteq \hat{\varphi}(p)} \kappa\setminus V^\alpha_t $$ by monotonicity of $\langle V^\beta_s \colon s \in [\kappa]^{<\omega}\rangle$ $$\subseteq \kappa \setminus V^\alpha_{\hat{\varphi}(p)} \subseteq \kappa \setminus D^\alpha_{\hat{\varphi}(p)}.$$ Thus $|\bigcap_{\delta\in p}W^\alpha_\delta \cap D^\alpha_{\hat{\varphi}(p)}| < \kappa$. The proof is complete. \end{proof}
\begin{lemma} Let $\kappa$ be a regular cardinal and let $\hat{\varphi}$ be a $\kappa$-shrinking function. If $\{\mathcal{F}_{\xi, \zeta} \colon \xi, \zeta < \kappa\}$ is a stratified set of ultrafilters with property (P), then each $\mathcal{F}_{\xi, \zeta}, (\xi, \zeta < \kappa)$ is a Simon point. \end{lemma} \begin{proof} The set $\{\mathcal{F}_{\xi, \zeta} \colon \xi, \zeta < \kappa\}$ is stratified, hence each $\mathcal{F}_{\xi, \zeta}$ is a non-minimal ultrafilter. Let $\{X_\xi \colon \xi < \kappa\}$ be a $\kappa$-discrete set of ultrafilters such that $\mathcal{F}_{\xi, \zeta} \in \overline{X}_\xi \setminus X_\xi$. Let $$Y = \{\mathcal{G}_t \colon t < [\kappa]^{< \omega}\}$$ be a $\kappa$-discrete set of minimal ultrafilters in Rudin-Frolik order. Fix any partition $$\{D_t \colon t < [\kappa]^{< \omega}\}$$ of $\kappa$ such that $D_t \in \mathcal{G}_t$. Let $$W = \{\mathcal{F}_{\xi, \zeta} \colon \exists_{\alpha< \kappa} D_\alpha \in \mathcal{F}_{\xi, \zeta}\}.$$ To complete the proof, it is enough to show that $\tilde{W} \cap \overline{Y} = \emptyset.$ Obviously $W(0)\cap \overline{Y} = \emptyset,$ where $W(0) = W$.

Now, we will proceed by induction on $\gamma < 2^\kappa$. Assume that for some $\beta < \gamma$ we have $W(\beta)\cap \overline{Y} = \emptyset$ and there exists $\mathcal{F}_{\xi, \zeta}\in W(\beta+1)\cap \overline{Y}.$ Thus, there exists $D_t \in \mathcal{F}_{\xi, \zeta}$ and $\eta < \xi$ such that $$\{\mathcal{F}_{\eta, \nu} \colon D_t \in \mathcal{F}_{\eta, \nu}\} \subseteq W(\beta) \textrm{ for } \nu < \kappa.$$ Hence $\mathcal{F}_{\xi, \zeta} \in \overline{W(\alpha) \cap X_\eta}$ for $\eta>\xi$. Then $$\overline{W(\beta) \cap X_\eta \cap \overline{Y}} \not = \emptyset.$$ A contradiction by Fact 1 and the definition of non-minimal ultrafilters. \\ Thus, if $\mathcal{F}_{\xi, \zeta}\in \tilde{W}$, then $\mathcal{F}_{\xi, \zeta} \not \in \overline{Y}$.

Assume that $\mathcal{F}_{\xi, \zeta} \not \in \tilde{W}$. Then the property $(P)$ is fulfilled, (by Lemma~1). Let $\{W_\delta \colon \delta < \kappa^+\} \in \mathcal{F}_{\xi, \zeta}$ be such that $$|\bigcap_{\delta\in p}W_\delta \cap D_{\hat{\varphi}(p)}| < \kappa$$ for any $p \in [\kappa^+]^{<\omega}$. Then, there are $\delta \in p$ such that $W_\delta \not \in \mathcal{G}_{\hat{\varphi}(p)}.$ Thus $\mathcal{F}_{\xi, \zeta} \not \in \overline{Y}$. \end{proof} 
\begin{theorem} Let $\kappa$ be a regular cardinal and let $\hat{\varphi}$ be a $\kappa$-shrinking function. There exists a Simon point in $\beta \kappa \setminus \kappa$. \end{theorem} \begin{proof} The theorem implies directly from Lemma 1 and Lemma 2. \end{proof} \begin{corollary} Let $\kappa$ be a regular cardinal. There are $\kappa$ many $\kappa$-discrete sets $X_\alpha, \alpha < \kappa$ of Simon points such that $X_\beta = \overline{X}_\alpha \setminus X_\alpha$, for $\beta< \alpha< \kappa$. \end{corollary}

\noindent
\textbf{Open problem} Is there a Simon point in $\beta\kappa \setminus \kappa$ for $\kappa$-singular?

		\noindent
		{\sc Joanna Jureczko}
		\\
		Wroc\l{}aw University of Science and Technology, Wroc\l{}aw, Poland
		\\
		{\sl e-mail: joanna.jureczko@pwr.edu.pl}
\end{document}